\documentclass[a4paper, 12pt]{article}
\usepackage{latexsym}
\usepackage{amsmath}
\usepackage{graphics}
\usepackage{amsthm}
\usepackage{amssymb}
\usepackage[all]{xy}
\newfont{\Fr}{eufm10}
\newfont{\Sc}{eusm10}
\newfont{\Bb}{msbm10}
\newfont{\Am}{msam10}
\newfont{\am}{msam7}
\numberwithin{equation}{section}
\newtheorem{theorem}{Theorem}[section]

\newtheorem{lemma}[theorem]{Lemma}
\newtheorem{corollary}[theorem]{Corollary}

\newtheorem{ftheorem}{Theorem}{\bf}{\it}
{\bf}{\it}
\theoremstyle{definition}
\newtheorem{definition}[theorem]{Definition}

{\bf}{\rm}

\theoremstyle{remark}

\newtheorem{remark}[theorem]{Remark}
\newtheorem{definition and corollary}[theorem]{Definition and Corollary}

{\it}{\rm}
\newcommand{\Ext}{\mbox{\rm Ext}}

\newcommand{\Lie}{\mbox{\rm Lie }}

\newcommand{\MID}{\! \! \mid}

\newcommand{\h}{\mathfrac{\h}}

\newcommand{\Q}{\mathbb{Q}}

\newcommand{\Z}{\mathbb{Z}}

\title{On the combinatorics of unramified admissible modules}
\author{Syu \textsc{Kato}
          \footnote{Graduate School of Mathematical Sciences, University of Tokyo, 3-8-1 Meguro Komaba 153-8914, Japan.} \footnote{Supported by JSPS Research Fellowship for Young Scientists}}
\begin{document}
\maketitle

\begin{abstract}
We construct a certain topological algebra $\Ext ^{\sharp} _{G ^{\vee}} X ( \chi )$ from a Deligne-Langlands parameter space $X ( \chi )$ attached to the group of rational points of a connected split reductive algebraic group $G$ over a non-Archimedean local field $\mathbb K$. Then we prove the equivalence between the category of continuous modules of $\Ext ^{\sharp} _{G ^{\vee}} X ( \chi )$ and the category of unramified admissible modules of $G ( \mathbb K )$ with a generalized infinitesimal character corresponding to $\chi$. This is an analogue of Soergel's conjecture which concerns the real reductive setting.
\end{abstract}

\section*{Introduction}
Let $G (\mathbb K)$ be the group of rational points of a connected split reductive algebraic group over a non-Archimedean local field $\mathbb K$, let $I$ be an Iwahori subgroup of $G (\mathbb K)$, and let $\mathbb H$ be the Hecke algebra over $\mathbb C$ associated to $( G (\mathbb K), I )$. By the Borel-Casselman theorem, the category of unramified admissible modules $\mathcal M ( G (\mathbb K) ) ^{ur}$ of $G (\mathbb K)$ is equivalent to the category of finite dimensional $\mathbb H$-modules.
Let $G ^{\vee}$ be the connected reductive algebraic group over $\mathbb C$ which is dual to $G$ in the sense of Langlands. Then, simple objects of $\mathcal M ( G (\mathbb K) ) ^{ur}$ (and hence $\mathbb H$-mod) are roughly parametrized by $G ^{\vee}$-orbits of the set: 
$$X =\left\{(s,u)\in G ^{\vee} \times G ^{\vee} \mid \begin{matrix} s\;{\rm is\;semisimple}, \\u\;{\rm is\;unipotent},\\s us^{-1}=u^{q}\end{matrix}\right\},$$
where $q$ is the order of the residue field of $\mathbb K$ (so-called the Deligne-Langlands conjecture proved by Kazhdan-Lusztig \cite{KL2} and Ginzburg). Here specifying an infinitesimal character $\chi$ of $\mathbb H$ ($\cong$ $G ^{\vee}$-conjugacy class of $s$) amounts to choose a Deligne-Langlands parameter space $X ( \chi ) \subset X$. (For the original formulation, see Borel \cite{Bo}.) Hence, the study of simple objects of $\mathcal M ( G (\mathbb K) ) ^{ur}$ breaks into the pieces $X ( \chi )$.\\
In real reductive setting, the Deligne-Langlands parameter space was modified by Adams-Barbasch-Vogan \cite{ABV} in order to incorporate Vogan's character duality into the Langlands philosophy. Inspired by their work, Soergel \cite[2.4]{S3} has introduced the {\it geometric extension algebra} $\Ext ^{\sharp} _{G ^{\vee}} X ( \chi )$. Let $Y$ be a $G ^{\vee}$-orbit in $X ( \chi )$. For each irreducible $G ^{\vee}$-local system $\tau$ on $Y$, we have its minimal extension $IC ( Y, \tau )$ along $X ( \chi )$. We define $\mathcal P := \bigoplus IC ( Y, \tau )$, where the summation is taken in all $G ^{\vee}$-orbits $Y$ and all $G ^{\vee}$-equivariant irreducible local system $\tau$ on $Y$. We define
$$\Ext ^{\sharp} _{G ^{\vee}} X ( \chi ) := \prod _{n \ge 0} \Ext ^{n} _{D ^+ _{G ^{\vee}} ( X ( \chi ) )} ( \mathcal P, \mathcal P ).$$
Here we regard it as a topological algebra by using degree $\ge i$-parts as its open base. Let $\mathcal M ( G (\mathbb K) ) ^{ur} _{\chi}$ be the fullsubcategory of $\mathcal M ( G (\mathbb K) ) ^{ur}$ consisting of modules with the generalized infinitesimal character $\chi$. We prove an analogue of the conjecture of Soergel (cf. \cite{S3} Basic conjecture 1), which claims that $X ( \chi )$ carries the complete information of $\mathcal M ( G (\mathbb K) ) ^{ur} _{\chi}$:
\begin{ftheorem}[(Unramified $p$-adic analogue of Soergel's conjecture)]\label{SC}
We have the following category equivalence:
$$\mathcal M ( G (\mathbb K) ) ^{ur} _{\chi} \stackrel{\cong}{\longrightarrow} \Ext ^{\sharp} _{G ^{\vee}} X ( \chi ) \text{\rm -Nil}.$$
Here $A ^{\sharp}\text{\rm -Nil}$ denote the category of finite dimensional continuous $A ^{\sharp}$-modules of a topological ring $A ^{\sharp}$.
\end{ftheorem}
Let $s$ be a semisimple element of $G ^{\vee}$ which is contained in the image of the first projection of $X ( \chi )$. We put $G ^s := Z _{G ^{\vee}} ( s )$. Let $\mathcal N$ be the nilpotent cone of $G ^{\vee}$. We have an action of $G ^{\vee} \times \mathbb G _m$ on $\mathcal N$, where the second factor acts as the dilation. Let $\mathcal N ^s \subset \mathcal N$ denote the fixed point set of the $( s, q ^{- 1} )$-action. Then, we have $X ( \chi ) \cong G ^{\vee} \times _{G ^s} \mathcal N ^s$. Let $\pi ^s : \widetilde{\mathcal N} ^s \rightarrow \mathcal N ^s$ be the Springer resolution restricted to the fixed point set of the $( s, q ^{- 1} )$-action. In this setting, Theorem \ref{SC} is a corollary of the following due to induction equivalence. Note that as $q$ is not a root of unity, all $IC ( Y, \tau )$ appear in the pushforward $\pi _* ^s \mathbb C _{\widetilde{\mathcal N} ^s}$:
\begin{ftheorem}[(= Corollary \ref{equiv})]\label{intromain}
We have the following isomorphism of topological $\mathbb C$-algebras:
$$\mathbb H _{\chi} \stackrel{\cong}{\longrightarrow} \Ext ^{\sharp} _{G ^{s}} \left( \pi ^s _* \mathbb C _{\widetilde{\mathcal N} ^s}, \pi ^s _* \mathbb C _{\widetilde{\mathcal N} ^s} \right).$$
Here $\mathbb H _{\chi}$ is the topological completion of $\mathbb H$ at $\chi$.
\end{ftheorem}
This is the main result of this paper. If we drop the $G ^s$-equivariant structure, then we recover Ginzburg's description (cf. \cite[8.1.5]{CG}).

Roughly speaking, the proof of Theorem \ref{intromain} consists of three ingredients. One is to use some etale ring extensions to divide the study of both completed algebras into two parts which reflect the action of $\left< s \right>$. Another is to pursue a chain of geometric isomorphisms which is a mixture of \cite{CG} and \cite{Gi}. The third is to consider $q$ as an indeterminant and then to check that we can specialize to a particular value.

The counterpart of Vogan's character duality in $p$-adic setting is the Zelevinsky conjecture i.e. $p$-adic analogue of the Kazhdan-Lusztig multiplicity formula, which was settled by Ginzburg (cf. \cite[\S 8]{CG}).

{\bf Warning:} For the sake of simplicity, our exposition below totally forgets the original group $G$ in this introduction. In particular, the group $G$ in the later part of this paper was abbreviated as $G ^{\vee}$ above! Moreover, we denote $\mathbb R ^{\bullet} F$ (the right derived functor of $F$) by $F$.

\section{Notations}\label{not}
\baselineskip=6pt
\baselineskip=0.5\baselineskip
{\small
\begin{itemize}
\item $G$, $B$, and $T$ : A split simply-connected algebraic group over $\Q$, its Borel subgroup, and its maximal torus;
\item $( W, S )$ : the Coxeter group attached to $G, B, T$ ($W$ is the Weyl group of $(G, T)$ and $S$ is the set of simple reflections coming from $T \subset B \subset G$.);
\item $\ell : W \rightarrow \Z$ : the length function of $W$ with respect to $S$;
\item $\mathfrak l := \Lie L$ : the Lie algebra of an algebraic group $L$;
\item $R ( L )$ : the representation ring of $L$;
\item $\mathcal N$ : the set of (ad-)nilpotent elements of $\mathfrak g$ (the nilpotent cone);
\item $\widetilde{\mathcal N} := \{ ( n, g B ) \in \mathcal N \times G / B ; g n g ^{- 1} \in \Lie B \} \stackrel{\pi}{\rightarrow} \mathcal N$ (the Springer resolution);
\item $\mathcal B := G / B$ and $Z := \widetilde{\mathcal N} \times _{\mathcal N} \widetilde{\mathcal N}$ : the flag variety and the Steinberg variety;
\item $a := (s, q)$ : a fixed element of $G \times \mathbb G _m$, where $s = e ^{\lambda}$ ($\lambda \in \mathfrak t$) is semisimple;
\item $A$ : the Zariski closure of the group generated by $a$ in $G \times \mathbb G _m$;
\item $G ^s$ and $W ^s$ : the commutator group of $s$ in $G$ and its Weyl group regarded as a subgroup of $W$ via $N _{G ^s} ( T ) \subset N _G ( T )$;
\item $X ^a$ : the set of $a$-fixed points of a $G \times \mathbb G _m$-variety $X$.
\item $H ^{\bullet} \subset H ^{\sharp}$ : the total cohomology algebras $\oplus _n H ^{n} \subset \Pi _n H ^{n}$.
\item $\Ext ^{\bullet} _G$ : the (bi-)functor $\Ext ^{\bullet}$ in $G$-equivariant derived categories;
\item $R ^{\sharp} \text{\rm -Nil}$ : the category of finite dimensional $R ^{\bullet}$-modules which are annihilated by $R ^{\ge n}$ for $n >> 0$ ($R ^{\bullet}$ is a graded algebra). ($=$ the category of finite dimensional continuous $R ^{\sharp}$-modules.)
\end{itemize}}
\normalbaselines
\section{Statement of the results}

\begin{definition}[(Iwahori-Hecke algebras)]
\begin{enumerate}
\item A finite Hecke algebra $H _W$ is a free $\Z [ v, v ^{- 1} ]$-module with basis $\{ T _w \} _{w \in W}$ such that the following multiplication rule holds:
\begin{itemize}
\item $( T _{s _{\alpha}} + 1 ) ( T _{s _{\alpha}} - v ) = 0$ if $s _{\alpha} \in S$;
\item $T _y \cdot T _w = T _{yw}$ if $\ell ( y ) + \ell ( w ) = \ell ( y w )$.
\end{itemize}
\item An affine Hecke algebra $\mathbb H$ is a free $\Z [ v, v ^{- 1} ]$-module with basis $\{ e ^{\lambda} \cdot T _w ; w \in W, e ^{\lambda} \in R ( T ) \}$ such that:
\begin{itemize}
\item $\{T _w\}$ span a subalgebra of $\mathbb H$ isomorphic to $H _W$;
\item $\{e ^{\lambda} \}$ span a $\Z [ v, v ^{- 1}]$-subalgebra of $\mathbb H$ isomorphic to $R ( T ) [v, v ^{- 1}]$;
\item For each $s _{\alpha} \in S$ and $\lambda \in X ^* ( T )$, we have
$$e ^{s _{\alpha} ( \lambda )} T _{s _{\alpha}} - T _{s _{\alpha}} e ^{\lambda} = ( 1 - v ) \frac{e ^{\lambda} - e ^{s _{\alpha} ( \lambda )}}{1 - e ^{- \alpha}}.$$
\end{itemize}
\end{enumerate}
\end{definition}

\begin{theorem}[(Bernstein)]
The center of $\mathbb H$ is equal to $R ( T ) ^W [ v, v ^{- 1} ]$.
\end{theorem}

Let $ev _s$ be the evaluation map $R ( G ) \rightarrow \mathbb C$ at $s$. Let $\widehat{R ( G )} _s$ be the completion of $\mathbb C \otimes _{\mathbb Z} R ( G )$ with respect to (the maximal ideal) $\ker ev _{s}$. By using the isomorphism $R ( G ) \cong R ( T ) ^W$, we define the completed Hecke algebra (with respect to $a$) as
$$\widehat{\mathbb H} _a := \mathbb H \MID _{v = q} \otimes _{R ( G )} \widehat{R ( G )} _s.$$
For each $R ( G )$-algebra $\mathcal R$, the ring $\widehat{R ( G )} _s \otimes _{R ( G )} \mathcal R$ admits a natural topological structure induced from the topology of $\widehat{R ( G )} _s$. In particular, $\widehat{\mathbb H} _a$ is a toplogical ring. Consider the restriction $\pi ^a : \widetilde{\mathcal N} ^a \rightarrow \mathcal N ^a$ of $\pi$ to $a$-fixed points. $\pi ^a$ is a $G ^s \times \mathbb G _m$-equivariant morphism. We regard
$$\Ext _{G ^s} ^{\sharp} \left( \pi ^a _* \mathbb C _{\widetilde{\mathcal N} ^a}, \pi ^a _* \mathbb C _{\widetilde{\mathcal N} ^a} \right) := \prod _{n \ge 0} \Ext _{D ^+ _{G ^s} ( \mathcal N ^s )} ^{n} \left( \pi ^a _* \mathbb C _{\widetilde{\mathcal N} ^a}, \pi ^a _* \mathbb C _{\widetilde{\mathcal N} ^a} \right)$$
as an algebra via the Yoneda product (cf. \S \ref{edf}). By choosing
$$U _i := \prod _{n \ge i} \Ext _{G ^s} ^{n} \left( \pi ^a _* \mathbb C _{\widetilde{\mathcal N} ^a}, \pi ^a _* \mathbb C _{\widetilde{\mathcal N} ^a} \right) \text{ for all } i \in \mathbb N$$
as a fundamental open neighborhood system, we equip $\Ext _{G ^s} ^{\sharp} \left( \pi ^a _* \mathbb C _{\widetilde{\mathcal N} ^a}, \pi ^a _* \mathbb C _{\widetilde{\mathcal N} ^a} \right)$ a topological structure. Then, we prove:
\begin{theorem}\label{main}
We have the following isomorphism of topological $\mathbb C$-algebras:
$$\widehat{\mathbb H} _a \cong \Ext _{G ^s} ^{\sharp} \left( \pi _* ^a \mathbb C _{\widetilde{\mathcal N} ^a}, \pi _* ^a \mathbb C _{\widetilde{\mathcal N} ^a} \right).$$

\end{theorem}
\begin{remark}
Unfortunately, the construction of our isomorphism is not canonical. If we forget $G ^s$-equivariant structure as in (\ref{forget}), then we get a natural isomorphism between the specialized Hecke algebra and the $\Ext$-algebra due to Ginzburg \cite{CG} 8.1.5.
\end{remark}

Let $\widehat{\mathbb H} _a\text{-Nil}$ be the category of finite dimensional continuous $\widehat{\mathbb H} _a$-modules.

\begin{corollary}[(Unramified $p$-adic analogue of Soergel's conjecture cf. \cite{S3})]\label{equiv}
We have the following category equivalence:
$$\widehat{\mathbb H} _a\text{\rm -Nil} \stackrel{\cong}{\longrightarrow} \Ext _{G ^s} ^{\sharp} \left( \pi _* ^a \mathbb C _{\widetilde{\mathcal N} ^a}, \pi _* ^a \mathbb C _{\widetilde{\mathcal N} ^a} \right) \text{\rm -Nil}.$$
\end{corollary}

\section{Equivariant derived categories and related constructions}\label{edf}
In this section, we review some constructions which are needed in the sequel. We included this section since they are scattering in the literature. We basically follow the setting of Bernstein-Lunts \cite{BL} with the necessary supplements. Let $L$ be a reductive linear algebraic group acting on a smooth variety $M$. There exists a sequence of $L$-varieties
$$\emptyset \neq E L _0 \hookrightarrow E L _1 \hookrightarrow E L _2 \hookrightarrow \cdots$$
such that 1) $L$ acts on each $E L _i$ freely, 2) $E L _i$ is connected and locally connected, and 3) $H ^j ( E L _i, R ) = 0$ for every $0 < j < i$ and any ring $R$. We call the limit $E L := \varinjlim E L _i$ a classifying bundle of $L$. Then, we have
$$M \stackrel{p}{\longleftarrow} E L \times M \stackrel{q}{\longrightarrow} E L \times _L M =: M _L,$$
where $q$ is a map obtained by the diagonal quotient by $L$. We define the $L$-equivariant derived category as
$$D ^+ _L ( M ) := \{ ( F, \underline{F}, \alpha _M ) \in D ^+ ( M ) \times D ^+ ( M _L ) \times \mathcal{H}om _{D ^+ ( E L \times M )} ; \alpha _M : q ^* \underline{F} \cong p ^* F \},$$
where $D ^+ ( M )$ stands for the derived category of the category of $\mathbb C _M$-modules bounded from below. We have a forgetful functor $For : D ^+ _L ( M ) \rightarrow D ^+ ( M )$ defined as the first projection. We define the category of $L$-equivariant perverse sheaves $\mathrm{Perv} ^L M$ as the pullback of $\mathrm{Perv} M$ via $For$. (Here the perversity is middle. cf. \cite{KS} \S 10.2) We define $BL := E L / L$, which is equal to $M _L$ when $M = pt$. For each $L$-equivariant locally constant sheaf $F$, we define its $L$-equivariant cohomology as
$$H ^{\bullet} _L ( M, F ) := H ^{\bullet} ( M _L, \underline{F} ),$$
where $\underline{F}$ is a locally constant sheaf on $M _L$ such that $p ^* F \cong q ^* \underline{F}$. Here the fibration $M \hookrightarrow M _L \rightarrow B L$ induces the following spectral sequence
\begin{eqnarray}
E _2 ^{r, t} := H ^{t} ( B L, \mathbb C _{B L} ) \otimes H ^r ( M, F ) \Rightarrow H ^{t + r} _L ( M, F ).\label{forget}
\end{eqnarray}
In particular, $H ^{\bullet} _L ( M, F )$ is a free $H ^{\bullet} _L ( pt )$-module if $H ^\bullet ( M, F )$ is concentrated in odd or even degree. Hence, we have
$$\mathbb C \otimes _{H ^{\bullet} _L ( pt )} H ^{\bullet} _L ( M, F ) \cong H ^{\bullet} ( M, F )$$
in this case.

Let $\mathcal E$ be a $L$-equivariant locally free sheaf over $M$. Since $q$ is a smooth morphism (of finite relative dimension), there is a $fpqc$-descent locally free sheaf $\underline{\mathcal E}$ of $p ^* \mathcal E$ along $q$ (cf. Bosch, L\"utkebohmert, and Raynaud \cite[Chap. 6 Theorem 4]{BLR}). Let $\mathbf{Vect} ^L ( EL \times M )$ be the category defined as:
\begin{itemize}
\item {\bf (Objects)} Data $\mathcal F := \{ \mathcal F _i, \phi _{i, j} \} _{i \ge j}$ such that 1) $\mathcal F _i$ is a $L$-equivariant algebraic vector bundle on $E L _i \times M$, 2) $\phi _{i, j} : \mathcal F _i \otimes \mathcal O _{E L _j \times M} \stackrel{\cong}{\longrightarrow} \mathcal F _j$ is a $L$-equivariant isomorphism, and 3) $\phi _{j, k} \circ \phi _{i, j} = \phi _{i, k}$;
\item {\bf (Morphisms)} A morphism
$$f : \{ \mathcal F _i, \phi _{i, j} \} _{i \ge j} \rightarrow \{ \mathcal F _i ^{\prime}, \phi _{i, j} ^{\prime} \} _{i \ge j}$$
is a family of $L$-equivariant morphisms $f _i : \mathcal F _i \rightarrow \mathcal F _i ^{\prime}$ as coherent sheaves such that $f _j \circ \phi _{i , j} = \phi _{i, j} ^{\prime} \circ f _j$.
\end{itemize}
We define the $L$-equivariant $K$-group of an ind-scheme $E L \times M$ as the Grothendieck group of the category $\mathbf{Vect} ^L ( EL \times M )$. Similarly, we define the $K$-group of an ind-scheme $M _L$ by replacing $E L _i \times M$ with $E L _i \times _L M$ and drop the $L$-equivariance condition from every object and morphism. Then, we have a map $\mathrm{sw} : K ^L ( M ) \rightarrow K ^L ( EL \times M ) \cong K ( M _L )$. Once we compose $\mathrm{sw}$, every construction of characteristic classes, the Mayer-Vietris exact sequence, or the Riemann-Roch theorem fall into the corresponding statements for $M _L$. For example, from the Chern character map $ch ^* : K ( M _L ) \longrightarrow H ^{\sharp} ( M _L ),$
we define the $L$-equivariant Chern character map as
\begin{eqnarray}
ch ^* _L := ch ^* \circ \mathrm{sw}.\label{e-chern}
\end{eqnarray}

\section{Convolution algebras}\label{convalg}
We retain the setting of the previous section. In this section, we review the convolution construction. We included this section since they are unfortunately not explicitly presented in \cite{CG} or \cite{Gi}. Therefore, most of the materials are borrowed from \cite{CG} with the necessary modifications. Let $a \in L$ be a central element. Let $A$ denote the Zariski closure of $\left< a \right>$ in $L$. We denote the natural inclusion $M ^a \hookrightarrow M$ by $i$. Then we define
$$\lambda _a := i ^* i _* 1 \in K ^L ( M ^a ).$$
Let $R ( L ) _a$ be the localization of $R ( L )$ at $a$. We put $K ^L ( M ) _a := \mathbb C \otimes _{\mathbb Z} R ( L ) _a \otimes _{R ( G )} K ^L ( M )$. By Thomason's argument \cite[\S 6]{Th}, we define the map $\mathrm{res} _a$ as $( \lambda _a ^{- 1} \cdot ) \circ i ^* : K ^L ( M ) _a \rightarrow K ^L ( M ^a ) _a$. We have $res _a \circ i _* = id$.
\begin{theorem}[(Thomason \cite{Th} Corollary 6.7)]\label{comm}
Assume that $L$ is abelian. Let $f : M \rightarrow N$ be a $L$-equivariant proper morphism of smooth $L$-varieties. Then, the following diagram is commutative:
$$
\xymatrix{
K ^L ( M ) _a \ar[r]^{f _*} \ar[d]^{res _a} & K ^L ( N ) _a\ar[d]^{res _a}\\
K ^L ( M ^a ) _a \ar[r]^{f _*} & K ^L ( N ^a ) _a\\
}.
$$
\end{theorem}

We denote $M \times M$ by $\mathbb M$. Consider an $A$-stable closed subvariety $Z$ of $\mathbb M$. We have the following diagram. 
$$
\xymatrix{
Z ^a \ar@{^{(}->}[r]^{i} \ar@{^{(}->}[d] & Z \ar@{^{(}->}[d]\\
M ^a \times M ^a \ar@{^{(}->}[r] & M \times M \\
}.
$$
We define $r _a$ as follows:
$$r _a : K ^L ( Z ) _a \stackrel{i ^*}{\rightarrow} K ^L ( Z ^a ) _a \stackrel{1 \boxtimes \lambda _a ^{- 1}}{\rightarrow} K ^L ( Z ^a ) _a.$$
Let $p _{ij}$ ($1 \le i < j \le 3$) be the projection from $M ^{3}$ to its $(i,j)$-th factor. Assume that we have $p _{13} ( p ^{- 1} _{12} ( Z ) \cap p ^{- 1} _{23} ( Z ) ) \subset Z$ and $p _{13}$ is projective. Then we define the convolution product as follows:
$$\star : K ^L ( Z ) \otimes K ^L ( Z ) \ni (\mathcal F, \mathcal G) \mapsto \left( p _{13} \right) _* \left( p _{12} ^* \mathcal F \otimes _{\mathcal O _{M ^3}} p _{23} ^* \mathcal G \right) \in K ^L ( Z ).$$
\begin{theorem}[(cf. Theorem 5.11.10 of \cite{CG})]\label{map r}
Assume that $L$ is abelian. Then, the map $r _a$ commutes with convolution.
\end{theorem}
\begin{proof}
The proof is exactly the same as \cite{CG} Theorem 5.11.10 if we replace \cite{CG} Theorem 5.11.7 by Theorem \ref{comm}.
\end{proof}

We define the equivariant relative cohomology $H ^{\bullet} _L ( Z ^a )$ of $\mathbb M - Z ^a \subset \mathbb M$ by
$$H ^{\bullet} _L ( Z ^a \mid \mathbb M ) := H _{L} ^{\bullet} ( \mathbb M, \mathbb M - Z ^a ) = H ^{\bullet} \left( \text{Cone} [ \mathbb C _{\mathbb M} \rightarrow j _* \mathbb C _{\mathbb M - Z ^a} ] \text{ in } D ^+ _L ( \mathbb M ) \right),$$
where $j : \mathbb M - Z ^a \hookrightarrow \mathbb M$ is an inclusion. We denote the (relative) first projection map $EL \times _L \mathbb M \rightarrow EL \times _L (M \times pt)$ by $p _1 ^L$. Let $td ( p _1 ^L )$ be the relative Todd class with respect to $p _1 ^L$. Then we define the (equivariant) Riemann-Roch map $RR ^L$ as follows:
$$RR ^L : \mathbb C \otimes _{\mathbb Z} K ^L ( Z ^a ) \ni \mathcal F \mapsto td ( p _1 ^L ) \cup ch ^* _L \left( \mathcal F \right) \in H ^{\sharp} _L ( Z ^a \mid \mathbb M ).$$
We define the (non-equivariant) Riemann-Roch map by
$$RR : \mathbb C \otimes _{\mathbb Z} K ( Z ^a ) \ni \mathcal F \mapsto ( 1 \boxtimes td ( M ) ) \cup ch ^* \left( \mathcal F \right) \in H ^{\bullet} ( Z ^a \mid \mathbb M )$$
as in \cite{CG} p301.

\begin{corollary}\label{free coh}
Assume that $RR : \mathbb C \otimes _{\Z} K ( Z ^a ) \rightarrow H ^{\bullet} ( Z ^a \mid \mathbb M )$ is an isomorphism. Then, $H ^{\bullet} ( Z ^a \mid \mathbb M)$ is concentrated in the even degree.
\end{corollary}
\begin{proof}
The assertion follows from the fact that every algebraic cycle gives rise to an even-dimensional cocycle.
\end{proof}

We have a convolution algebra structure in $H ^{\bullet} _L ( Z ^a \mid \mathbb M )$ as in \cite[2.7]{CG} when $p _{13} : ( p _{12} ^{- 1} ( Z ^a ) \cap p _{23} ^{- 1} ( Z ^a ) ) \rightarrow Z ^a$ is projective.

\begin{theorem}[(cf. \cite{CG} Theorem 5.11.11)]\label{RR is alg hom}
$RR ^L$ commutes with convolution.
\end{theorem}
\begin{proof}
We have $sw \circ p _{ij} ^* = ( p _{ij} ^L ) ^* \circ sw$. Noting that the pullback $K ^L ( Z ^a ) \rightarrow K ^L ( EL \times Z ^a ) \cong K ( Z ^a _L )$ is an algebra morphism, the proof of \cite{CG} Theorem 5.11.11 works word by word.
\end{proof}

\begin{theorem}[(cf. \cite{CG} Theorem 8.6.7)]\label{int}
In the same notation as in \S 1, we have
$$H ^{\bullet} _{G ^s \times \mathbb G _m} ( Z ^a \mid \mathbb M ) \cong \Ext ^{\bullet} _{G ^s \times \mathbb G _m} \left( \pi ^a _* \mathbb C _{\widetilde{\mathcal N} ^a}, \pi ^a _* \mathbb C _{\widetilde{\mathcal N} ^a} \right)$$
as algebras. Here the product of the LHS is the convolution product and the product of the RHS is the Yoneda product.
\end{theorem}
\begin{proof}
This is a relative version of \cite[Theorem 8.6.7]{CG}. This theorem is a special version of \cite[Theorem 8.6.35]{CG}, which uses 1) adjunction maps, 2) tautological embeddings, 3) the Verdier duality, 4) the spectral sequence associated to the composition of two direct images, 5) the base change theorem, 6) $f _! = f _*$ for a proper map $f$, 7) commutativity of diagram. Hence, the result follows since all of the above constructions have counterparts for fibrations over $B ( G ^s \times \mathbb G _m )$.
\end{proof}

\begin{remark}
As in the arguments in \cite{CG}, we always take the support of cycles into account. Hence, we can apply the above technique for singular varieties as well (by an inexplicit mention to the ambient spaces). Here we normally use $M := \widetilde{\mathcal N}$.
\end{remark}

\section{Proof of Theorem \ref{main}}
Our starting point is the following theorem:
\begin{theorem}[(Ginzburg \cite{G})]
$K ^{G \times \mathbb G _m} ( Z ) \cong \mathbb H$ via the convolution product.
\end{theorem}

\begin{lemma}\label{algebra induction}
We have $\phi : R ( G ^s ) \otimes _{R ( G )} K ^{G \times \mathbb G _m} ( Z ) \cong K ^{G ^s \times \mathbb G _m} ( Z )$ as $R ( G ^s ) [ v, v ^{-1} ]$-algebras.
\end{lemma}
\begin{proof}
Since $G$ is a semisimple simply connected group, $R ( T )$ is a free $R ( G )$-module by the Pittie-Steinberg theorem. By the K\"unneth theorem (\cite{CG} Theorem 5.6.1 (d) $\Rightarrow$ (a)), we have an isomorphism
$$K ^{G} ( Y ) \otimes _{R ( G )} K ^G ( \mathcal B ) \rightarrow K ^G ( Y \times \mathcal B )$$
for an arbitrary $G$-variety $Y$. Hence, we have a chain of isomorphisms
\begin{eqnarray*}
R ( G ^s ) \otimes _{R ( G )} K ^{G \times \mathbb G _m} ( E ) \cong R ( G ^s ) \otimes _{R ( G )} K ^{G \times \mathbb G _m} ( \mathcal B ) \qquad\\
\qquad \cong K ^{G \times \mathbb G _m} ( G / G ^s \times \mathcal B ) \cong K ^{G ^s \times \mathbb G _m} ( \mathcal B ) \cong K ^{G ^s \times \mathbb G _m} ( E )
\end{eqnarray*}
for each $G \times \mathbb G _m$-equivariant vector bundle $E$ over $\mathcal B$ by the Thom isomorphism and the induction argument. Therefore, the desired isomorphism as $R ( G ^s )$-modules follows from the localization sequence as in \cite[6.2.3]{CG}. Here $R ( G ) \subset K ^{G \times \mathbb G _m} ( Z )$ is contained in the center of $\mathbb H$. Hence, $\phi$ is a coefficient extension of algebras.
\end{proof}

\begin{lemma}\label{etaleness}
The ring extension $R ( G ) \subset R ( G ^s )$ is etale over $s \in T$.
\end{lemma}
\begin{proof}
The stabilizer of $s$ in $W$ is $W ^s$. Hence, we have an etale map $T / W \rightarrow T / W ^s$, which is the spectrum of the ring inclusion.
\end{proof}

By definition, we have $a = ( s, q ) \in G ^s \times \mathbb G _m$.

\begin{corollary}\label{f}
We have an isomorphism $\widehat{K ^{G \times \mathbb G _m} ( Z )} _a \stackrel{\cong}{\longrightarrow} \widehat{K ^{G ^s \times \mathbb G _m} ( Z )} _a$
as topological $\widehat{R ( G ^s \times \mathbb G _m )} _a$-algebras.
\end{corollary}
\begin{proof}
Combine Lemma \ref{algebra induction} and Lemma \ref{etaleness}.
\end{proof}

By Segal-Thomason's localization theorem \cite[5.3]{Th}, we have
$$i _* : K ^{T \times \mathbb G _m} ( Z ^a ) _a \cong K ^{T \times \mathbb G _m} ( Z ) _a.$$
Here the both sides are free $\mathbb C \otimes _{\mathbb Z} R ( T ) _a$-modules.
\begin{theorem}[(\cite{CG} Corollary 8.8.9)]
A connected component of $\widetilde{\mathcal N} ^a$ is a $G ^s \times \mathbb G _m$-equivariant vector bundle over the flag variety of $G ^s$.
\end{theorem}

\begin{corollary}\label{Z}
The subspace $Z ^a \subset Z$ is stratified by $G ^s \times \mathbb G _m$-equivariant vector bundles over $G ^s$-orbits of the two-fold product of flag varieties of $G ^s$.
\end{corollary}

\begin{lemma}
Under the above settings, we have
$$i _* : K ^{G ^s \times \mathbb G _m} ( Z ^a ) _a \cong K ^{G ^s \times \mathbb G _m} ( Z ) _a.$$
\end{lemma}

\begin{proof}
By counting the number of $T$-fixed points, we deduce that the stratification of Corollary \ref{Z} consists of $\left| W \right| ^2 / \left| W ^s \right|$ vector bundles over $G ^s$-orbits of the two-fold product of the flag variety of $G ^s$. By the Thom isomorphism and the localization sequence, $K ^{G ^s \times \mathbb G _m} ( Z ^a )$ is a free $R ( T )$-module of rank $\left| W \right| ^2 / \left| W ^s \right|$. We have a stratification $\mathcal S$ of the two-fold product of flag varieties of $G$ into $\left| W \right| ^2 / \left| W ^s \right| ^2$ vector bundles over a two-fold copies of the flag variety of $G ^s$ by collecting the attracting set of $s$. Consider its subdivision $\mathcal S ^+$ formed by the attracting set of the diagonal $G ^s$-orbits of $\mathcal B ^s \times \mathcal B ^s$. Then, $Z$ is stratified by vector bundles over strata in $\mathcal S ^+$. It follows that $K ^{G ^s \times \mathbb G _m} ( Z )$ is a free $R ( T )$-module of rank $\left| W \right| ^2 / \left| W ^s \right|$.
Consider the following commutative diagram:
$$
\xymatrix{
K ^{T \times \mathbb G _m} ( Z ) _a \ar[r] ^{res _a} & K ^{T \times \mathbb G _m} ( Z ^a ) _a \ar[r] ^{i _*} & K ^{T \times \mathbb G _m} ( Z ) _a\\
R ( T ) _a \otimes _{R ( G ^s ) _a} K ^{G ^s \times \mathbb G _m} ( Z ) _a \ar[r] ^{res _a} \ar[u] ^{\cong} & R ( T ) _a \otimes _{R ( G ^s ) _a} K ^{G ^s \times \mathbb G _m} ( Z ^a ) _a  \ar[u] \ar[r] ^{1 \times i _*} & R ( T ) _a \otimes _{R ( G ^s ) _a} K ^{G ^s \times \mathbb G _m} ( Z ) _a \ar[u] ^{\cong}
}
$$
Here vertical isomorphisms follow by Lemma \ref{algebra induction}. Since $i _* \circ res _a$ is the identity in the upper row, $R ( T ) _a \otimes _{R ( G ^s ) _a} K ^{G ^s \times \mathbb G _m} ( Z ^a ) _a$ contains $R ( T ) _a \otimes _{R ( G ^s ) _a} K ^{G ^s \times \mathbb G _m} ( Z ) _a$ as its direct summand. Therefore, the result follows since that the both are free modules of the same rank.
\end{proof}
Now Theorem \ref{map r} yields an algebra isomorphism $K ^{G ^s \times \mathbb G _m} ( Z ) _a \cong K ^{G ^s \times \mathbb G _m} ( Z ^a ) _a$ since they are subalgebras of $K ^{T \times \mathbb G _m} ( Z ) _a$.

Let $A ^c$ be a subgroup of $G ^s \times \mathbb G _m$ such that 1) $A A ^c = G ^s \times \mathbb G _m$ and 2) $A \cap A ^c$ is a finite group. This is possible because $A$ is central in $G ^s \times \mathbb G _m$.
\begin{lemma}\label{non-gen}
If $\mathcal N ^a \neq \{0\}$ and $q \neq \sqrt[l]{1}$, then we can choose $A ^c$ so that $A ^c \subset G ^s$.
\end{lemma}
\begin{proof}
Let $\mathfrak g _{\alpha} \subset \mathcal N ^a \subset \mathfrak g$ be a root space. Then, we have $e ^{\left< \alpha, \lambda \right>} = q$. This means $A T = T \times \mathbb G _m$. Hence, a Lie algebra calculation yields the result.
\end{proof}
We put $\underline{A} := ( G ^s \times \mathbb G _m ) / A ^c$. 
\begin{lemma}
The diagonal class of $\underline{A} \times \underline{A}$ is sitting in the image of
$$K ^{G ^s \times \mathbb G _m} ( \underline{A} ) \otimes _{R ( G ^s \times \mathbb G _m )} K ^{G ^s \times \mathbb G _m} ( \underline{A} ) \rightarrow K ^{G ^s \times \mathbb G _m} ( \underline{A} \times \underline{A} ).$$
\end{lemma}
\begin{proof}
Since $\underline{A}$ is a direct product of split tori and finite abelian group, we reduce the problem to the case $\underline{A}$ is either 1) a finite group or 2) $\mathbb G _m$ from \cite{CG} 5.6.1 (a) $\Leftrightarrow$ (b). Case 1) is trivial since the total space is discrete. Case 2) follows from the fact that the ideal $(x - y) \mathbb C [x, y] \subset \mathbb C [x, y]$ is in (diagonal) $\mathbb G _m$-equivariant linear equivalence with $x \mathbb C [x, y]$.
\end{proof}
It follows that
\begin{eqnarray*}
K ^{A ^c} ( Z ^a ) \cong K ^{G ^{s} \times \mathbb G _m} ( \underline{A} \times Z ^a ) \cong K ^{G ^{s} \times \mathbb G _m} ( \underline{A} ) \otimes _{R ( G ^s \times \mathbb G _m )} K ^{G ^{s} \times \mathbb G _m} ( Z ^a ) \\
 \cong R ( A ^c ) \otimes _{R ( G ^s \times \mathbb G _m )} K ^{G ^s \times \mathbb G _m} ( Z ^a )
\end{eqnarray*}
by induction and the K\"unneth theorem (cf. \cite{CG} Theorem 5.6.1 (b) $\Rightarrow$ (a)). Since $A$ acts on $Z ^a$ trivially, we obtain
$$K ^{A \times A ^c} ( Z ^a ) \cong R ( A \times A ^c ) \otimes _{R ( G ^s \times \mathbb G _m )} K ^{G ^s \times \mathbb G _m} ( Z ^a )$$
as $\mathbb C$-algebras.

\begin{lemma}\label{e}
The ring extension $R ( G ^s \times \mathbb G _m ) \subset R ( A \times A ^c )$ is etale along $a \in T \times \mathbb G _m$.
\end{lemma}
\begin{proof}
The product map $\theta : A \times A ^c \rightarrow G ^s \times \mathbb G _m$ is an etale map. Thus, its restriction to the maximal torus of $A \times A ^c$ which maps to $T \times \mathbb G _m$ is still etale. We denote it by $\theta ^{\prime}$. The Weyl group of $A ^c$ is $W ^s$, which is isomorphic to the Weyl group of $G ^s$. Hence, the $W ^s$-quotient of $\theta ^{\prime}$ is again etale. Therefore, the corresponding inclusion of coordinate rings yields the result.
\end{proof}
By abuse of notation, let us denote the point $(a, 1) \in A \times A ^c$ also by $a$. Then, we have an isomorphism
\begin{eqnarray}
\widehat{K ^{G ^s \times \mathbb G _m} ( Z ^a )} _a \stackrel{\cong}{\longrightarrow} \widehat{K ^{A \times A ^c} ( Z ^a )} _a.\label{group red}
\end{eqnarray}

\begin{lemma}\label{image of RR}
The equivariant Riemann-Roch map gives rise to an isomorphism of topological $\mathbb C$-algebras:
$$RR ^{A ^c} : \widehat{K ^{A ^c} ( Z ^a )} _1 \stackrel{\cong}{\longrightarrow} H ^{\sharp} _{A ^c} ( Z ^a \mid \mathbb M ).$$
\end{lemma}
\begin{proof}
We have the following commutative diagram:
$$
\xymatrix{
\mathbb C \otimes _{\mathbb Z} K ^{A ^c} ( Z ^a ) \ar[r]^{RR ^{A ^c}} \ar[d] & H ^{\sharp} _{A ^c} ( Z ^a \mid \mathbb M ) \ar[d]\\
\mathbb C \otimes _{\mathbb Z} K ( Z ^a ) \ar[r]^{RR} & H ^{\sharp} ( Z ^a \mid \mathbb M )
}
$$
Here the vertical arrows are surjective and the bottom arrow is an isomorphism (cf. \cite{CG} Proof of Proposition 8.1.5). Moreover, $K ^{A ^c} ( Z ^a )$ is a free $R ( A ^c )$-module by its construction and $H ^{\sharp} _{A ^c} ( Z ^a \mid \mathbb M )$ is a free $H ^{\sharp} _{A ^c} (pt)$-module by Corollary \ref{free coh}. The $A ^c$-equivariant Riemann-Roch map over a point is just expanding characters of $R ( A ^c )$ along $1$ since the relative Todd class is trivial in this case. Therefore, we have
$$RR ^{A ^c} : \widehat{K ^{A ^c} ( pt )} _1 \stackrel{\cong}{\longrightarrow} H ^{\sharp} _{A ^c} ( pt ).$$
Hence, we have the result from Theorem \ref{RR is alg hom}.
\end{proof}
We have a non-canonical isomorphism $\kappa : \widehat{R ( A )} _{a} \cong \mathbb C [[ \mathfrak a ]]$ since both sides are isomorphic to the completion of a polynomial algebra of dimension $\dim A$. Consider the diagonal action of $A \times A ^c$ on $E ( A \times A ^c ) \times E ( G ^s \times \mathbb G _m ) \times Z ^a$, which is homotopy equivalent to $E ( A \times A ^c ) \times Z ^a$. We can take $B ( A \cap A ^c ) = E ( A \times A ^c ) / A \cap A ^c$. Thus, we have a fibration
$$\vartheta : E ( A \times A ^c ) \times _{A \times A ^c} Z ^a \rightarrow E ( G ^s \times \mathbb G _m ) \times _{G ^s \times \mathbb G _m} Z ^a$$
with a fiber $B ( A \cap A ^c)$. Here we have $H ^{\bullet} ( B ( A \cap A ^c ) ) = \mathrm{Ext} ^{\bullet} _{\mathbb C [ A \cap A ^c]} \left( \mathbb C, \mathbb C \right)$. The trivial $( A \cap A ^c )$-module $\mathbb C$ is projective. Therefore, we conclude $H ^{\sharp} _{A \cap A ^c} (pt) \cong H ^{\sharp} ( B ( A \cap A ^c ) ) = \mathbb C$. Hence, the spectral sequence
$$E _2 := H ^{\bullet} _{A \cap A ^c} (pt) \otimes H ^{\bullet} _{G ^s \times \mathbb G _m} ( Z ^a \mid \mathbb M ) \Rightarrow H ^{\bullet} _{A \times A ^c} ( Z ^a \mid \mathbb M )$$
yields an isomorphism $\vartheta ^* : H ^{\sharp} _{G ^s \times \mathbb G _m} ( Z ^a \mid \mathbb M ) \cong H ^{\sharp} _{A \times A ^c} ( Z ^a \mid \mathbb M )$. $\vartheta ^*$ is an isomorphism of convolution algebras since $\vartheta$ is a $\mathbb C$-acyclic base extension. Therefore, we have a non-canonical chain of isomorphisms
\begin{eqnarray*}
\widehat{K ^{A \times A ^c} ( Z ^a )} _{a} \cong \widehat{R ( A )} _a \hat \otimes \widehat{K ^{A ^c} ( Z ^a )} _{1} \stackrel{\kappa \otimes RR ^{A ^c}}{\longrightarrow} \mathbb C [[ \mathfrak a ]] \hat \otimes H ^{\sharp} _{A ^c} ( Z ^a \mid \mathbb M ) \\
\cong H ^{\sharp} _{A \times A ^c} ( Z ^a \mid \mathbb M ) \cong H ^{\sharp} _{G ^s \times \mathbb G _m} ( Z ^a \mid \mathbb M ).
\end{eqnarray*}
Composing the all isomorphisms of topological $\mathbb C$-algebras given in Corollary \ref{f}, (\ref{group red}), and Theorem \ref{int}, we obtain
$$\widehat{K ^{G \times \mathbb G _m} ( Z )} _{a} \stackrel{\cong}{\rightarrow} \Ext ^{\sharp} _{G ^s \times \mathbb G _m} ( \pi ^a _* \mathbb C _{\widetilde{\mathcal N} ^a}, \pi ^a _* \mathbb C _{\widetilde{\mathcal N} ^a} ).$$
Now we want to specialize $v$ to a suitable value. If $A ^c$ is contained in $G ^s \times 1$, then substituting $v$ by $q$ at $K ^{A ^c} ( Z ^a )$ does nothing. We put $R := \left( R ( A ) \MID _{v = q} \right) _s \cap R ( A ^c ) _1$. By dimension counting, we have $\hat R \cong \mathbb C [[ \epsilon ]]$ or $\mathbb C$. In $R ( A )$ side, we have a non-canonical chain of isomorphisms
$$\kappa ^{\prime} : \widehat{\left( R ( A ) \MID _{v = q} \right)} _s \cong \widehat{R ( A / \mathbb G _m )} _1 \hat\otimes \hat{R} \cong \mathbb C [[ \mathfrak a / \Lie \mathbb G _m ]] \hat\otimes \hat{R},$$
where $\mathbb G _m \hookrightarrow A$ is an arbitrary embedding. Therefore, we have the following isomorphism
$$\widehat{\mathbb H} _a \cong \widehat{\left( R ( A ) \MID _{v = q} \right)} _s \hat\otimes _{\hat R} \widehat{K ^{A ^c} ( Z ^a )} _1 \stackrel{\kappa ^{\prime} \otimes RR ^{A ^c}}{\longrightarrow} \mathbb C [[ \mathfrak a / \Lie \mathbb G _m]] \hat\otimes H ^{\sharp} _{A ^c} ( Z ^a \mid \mathbb M ) \cong H ^{\sharp} _{G ^s} ( Z ^a \mid \mathbb M ).$$
Hence, Lemma \ref{non-gen} yields the result if $\mathcal N ^a \neq \{ 0 \}$ and $q \neq \sqrt[l]{1}$. If $q = \sqrt[l]{1}$, then "substituting" $v$ by $q$ at $K ^{A ^c} ( Z ^a )$ is the same as taking quotient with respect to the ideal $( v - 1 )$. Thus, it is compatible with the Riemann-Roch map. Since we have a non-canonical isomorphism $\widehat{R ( A )} _a \cong \mathbb C [[ \mathfrak a ]]$, the result follows in this case. If $\mathcal N ^a = \{ 0 \}$, then $1 \times \mathbb G _m$ acts on $Z ^a$ trivially. Hence, we have $K ^{G ^s \times \mathbb G _m} ( Z ^a ) \cong K ^{G ^s} ( Z ^a ) \otimes \mathbb Z [ v, v ^{- 1} ]$ as algebras. It follows that we can construct an isomorphism
$$\widehat{K ^{G ^s \times \mathbb G _m} ( Z ^a )} _a \cong \mathbb C [[ \epsilon ]] \hat{\otimes} \widehat{K ^{G ^s} ( Z ^a )} _{s}$$
in a non-canonical way. Similarly, we have
$$H ^{\sharp} _{G ^s \times \mathbb G _m} ( Z ^a \MID \mathbb M ) \cong \mathbb C [[ \epsilon ]] \hat{\otimes} H ^{\sharp} _{G ^s} ( Z ^a \MID \mathbb M )$$
as topological $\mathbb C$-algebras. Hence, we have
$$\widehat{( K ^{G ^s \times \mathbb G _m} ( Z ^a ) \MID _{v = q} )} _s \cong \widehat{K ^{G ^s} ( Z ^a )} _{s} \stackrel{\cong}{\longrightarrow} H ^{\sharp} _{G ^s} ( Z ^a \MID \mathbb M )$$
as topological $\mathbb C$-algebras by applying the similar arguments as above with $G ^s \times \mathbb G _m$ replaced by $G ^s$. Therefore, Theorem \ref{int} yields the result.

%{\bf Acknowledgement.} I am doing this work at Freiburg University.
{\bf Acknowledgement.} The author wants to thank Professor Wolfgang Soergel for letting me notice this problem and much discussions. The author wants to thank Professor Friedrich Knop for letting my stay possible, much discussions, and his tolerance. This work was done during my stay at Freiburg University in the fall semester 2003/2004. The author wants to thank the hospitality of the Institute and Freiburg city.

\end{document}